\documentclass[twoside,notitlepage,11pt]{article}

\usepackage{amssymb}
\usepackage[leqno]{amsmath}
\usepackage{amsfonts}
\usepackage{amsopn}
\usepackage{amstext}
\usepackage{amsthm}
\usepackage[all]{xy}

\providecommand{\cal}{\mathcal}
\renewcommand{\Bbb}{\mathbb}
\renewcommand{\frak}{\mathfrak}
\newenvironment{pf}{\begin{proof}}{\end{proof}}

\newcommand{\Aaa}{{\cal{A}}}
\newcommand{\Bee}{{\cal{B}}}
\newcommand{\Cee}{{\cal{C}}}

\newcommand{\Ef}{{\cal{F}}}

\newcommand{\Kay}{{\cal{K}}}
\newcommand{\El}{{\cal{L}}}
\newcommand{\Pee}{{\cal{P}}}

\newcommand{\Yu}{{\cal{U}}}

\newcommand{\Em}{{\mathcal{M}}}

\newcommand{\Nat}{{\Bbb{N}}}
\newcommand{\Qyu}{{\Bbb{Q}}}

\newcommand{\lam}{{\lambda}}
\newcommand{\al}{\alpha}

\renewcommand{\phi}{\varphi}
\renewcommand{\rho}{\varrho}

\newcommand{\rest}{\restriction}

\newcommand{\ntr}{n\in\omega}

\newcommand{\loe}{\leqslant}
\newcommand{\goe}{\geqslant}

\newcommand{\subs}{\subseteq}
\newcommand{\sups}{\supseteq}
\newcommand{\nnempty}{\ne\emptyset}

\newcommand{\ovr}{\overline}
\renewcommand{\iff}{\Longleftrightarrow}

\newcommand{\cl}{\operatorname{cl}}
\newcommand{\Int}{\operatorname{int}}
\newcommand{\w}{\operatorname{w}}

\newcommand{\id}{\operatorname{id}}

\newcommand{\liminv}{\varprojlim}
\newcommand{\lex}{\operatorname{lx}}

\newcommand{\poset}{{\Bbb{P}}}

\newcommand{\Ord}{{\frak{Ord}}}
\newcommand{\by}{/}

\newcommand{\setof}[2]{\{#1\colon #2\}}
\newcommand{\bigsetof}[2]{\Bigl\{#1\colon #2\Bigr\}}

\newcommand{\seqof}[2]{\langle #1\colon #2\rangle}
\newcommand{\sett}[2]{\{#1\}_{#2}}
\newcommand{\sn}[1]{\{#1\}} 
\newcommand{\dn}[2]{\{#1,#2\}} 
\newcommand{\pair}[2]{\langle #1, #2 \rangle} 
\newcommand{\map}[3]{#1\colon #2 \to #3} 
\newcommand{\img}[2]{#1[#2]} 
\newcommand{\inv}[2]{{#1}^{-1}[#2]} 

\newcommand{\dpower}[2]{[#1]^{#2}}

\newcommand{\fin}[1]{[#1]^{<\omega}}

\newcommand{\rloe}{\preceq}

\newcommand{\suppt}{\operatorname{suppt}}

\newcommand{\Sig}{\Sigma}

\providecommand{\nat}{\omega}

\newcommand{\R}{\ensuremath{\mathcal R}}

\newcommand{\RC}{\ensuremath{{\mathcal R \mathcal C}}}

\newcommand{\invsys}[5]{\langle {#1}_{#4},{#2}_{#4}^{#5},#3 \rangle}

\renewcommand{\S}{\mathbb S}

\newcommand{\rk}{\operatorname{rk}}

\newcommand{\cube}[1]{{[0,1]^{#1}}}

\newcommand{\closed}{\operatorname{Closed}}
\newcommand{\bbL}{{\mathbb L}}
\newcommand{\ray}{{\ensuremath{\operatorname{R}^{\rightarrow}}}}
\newcommand{\longint}{{\ensuremath{{}^{\leftarrow}\operatorname{I}^{\rightarrow}}}}

\newcommand{\ciag}[1]{\setof{{#1}_n}{\ntr}}

\newtheorem{tw}{Theorem}[section]
\newtheorem{wn}[tw]{Corollary}
\newtheorem{lm}[tw]{Lemma}
\newtheorem{prop}[tw]{Proposition}

\newtheorem{twkm}{Theorem}

\theoremstyle{definition}

\theoremstyle{remark}

\newtheorem{question}{Question}

\title{Compact spaces generated by retractions}
\author{
{\sc Wies{\l}aw Kubi\'s}
\thanks{Research supported by NATO Science Fellowship (March -- August 2004).}
\\
Instytut Matematyki, Akademia \'Swi\c etokrzyska\\
ul. \'Swi\c etokrzyska 15, 25-406 Kielce, Poland\\
E-mail: \texttt{wkubis@pu.kielce.pl}
}

\begin{document}
\maketitle

\begin{abstract} We study compact spaces which are obtained from metric compacta by iterating the operation of inverse limit of continuous sequences of retractions. This class, denoted by \R, has been introduced in \cite{BKT}. Allowing continuous images in the definition of class \R, one obtains a strictly larger class, which we denote by \RC. We show that every space in class \RC~is either Corson compact or else contains a copy of the ordinal segment $\omega_1+1$. This improves a result of Kalenda from \cite{Kalenda_segment}, where the same was proved for the class of continuous images of Valdivia compacta. We prove that spaces in class \R~do not contain {\em cutting P-points} (see the definition below), which provides a tool for finding spaces in $\RC\setminus\R$. Finally, we study linearly ordered spaces in class \RC. We prove that scattered linearly ordered compacta belong to $\RC$ and we characterize those ones which belong to \R. We show that there are only 5 types (up to order isomorphism) of connected linearly ordered spaces in class \R~and all of them are Valdivia compact. Finally, we find a universal pre-image for the class of all linearly ordered Valdivia compacta.


\end{abstract}

\section{Introduction}

Denote by \R~the smallest class of compact spaces containing all metric compacta and closed under limits of continuous inverse sequences whose bonding maps are retractions. This class has been introduced in \cite{BKT}, motivated by new results on locally uniformly convex renormings of Banach spaces and by new examples of compacta whose spaces of continuous functions have such a renorming. 
It has been proved in \cite{BKT} that $C(X)$ has an equivalent locally uniformly convex norm for every $X\in\R$. Actually, the argument given in \cite{BKT} shows that the same is true for spaces in a larger class, namely the smallest class that contains all metric compacta and which is closed both under continuous images and under limits of inverse sequences of retractions. We denote this class by \RC. We prove that countably tight spaces in class \RC~are Corson compact. More precisely, we show that if $X\in\RC$ is not Corson then the ordinal $\omega_1+1$ embeds into $X$. This generalizes the result of Kalenda from \cite{Kalenda_segment}, which deals with continuous images of Valdivia compacta.

Recall that a space $X$ is {\em Valdivia compact} \cite{AMN,DG} if for some $\kappa$ there exists an embedding $X\subs\cube\kappa$ such that $X = \cl (X\cap \Sigma(\kappa))$, where $\Sigma(\kappa)=\setof{x\in\cube \kappa}{|\suppt(x)|\loe\aleph_0}$ is the $\Sigma$-product of $\kappa$ copies of $\cube{}$ and $\suppt(x)=\setof{\al}{x(\al)\ne0}$. Let us mention that one of the important functional-analytic properties of Valdivia compacta is the existence of an equivalent locally uniformly convex norm on their spaces of continuous functions, see Chapter VII of \cite{DGZ}.
For other results on Valdivia compacta we refer to Kalenda's survey article \cite{Kalenda} and to the recent papers \cite{KM, KU}.

An important result from \cite{AMN} says that given a Valdivia compact $X$ suitably embedded in $\cube\kappa$, for every infinite set $S\subs \kappa$ there exists $T\sups S$ such that $|S|=|T|$ and the map $\map{r_T}X{\cube\kappa}$, defined by $r_T(x)(\al)=x(\al)$ for $\al\in T$ and $r_T(x)(\al)=0$ for $\al\in\kappa\setminus T$, is an internal retraction, i.e. $r_T(r_T(x))=r_T(x)\in X$ for every $x\in X$. This implies that every Valdivia compact space is the limit of a continuous inverse sequence of smaller Valdivia compacta whose bonding maps are retractions. In particular, Valdivia compacta belong to \R. It is an open question whether the class of Valdivia compacta is stable under retractions. We show that retracts of Valdivia compacta belong to class \R.

The study of classes \R~and \RC~can also be motivated by the following result from \cite{KM}:

\begin{twkm}[cf. {\cite[Corollary 4.3]{KM}}]\label{KM}
Valdivia compacta of weight $\loe\aleph_1$ are precisely those spaces which can be obtained as limits of continuous inverse sequences of metric compacta with right-invertible bonding maps.
\end{twkm}

Class \R~contains spaces which are not continuous images of Valdivia compacta. Perhaps the simplest example is the Alexandrov duplication of a countable dense set in the Cantor cube $2^{\aleph_1}$, which is a non-metrizable compactification of the natural numbers, see \cite[Example 4.6(b)]{KM}. Another, non-trivial, example is a compact linearly fibered ccc non-separable space $K$ constructed by Todor\v cevi\'c (see the proof of  Theorem 8.4 in \cite{To_chains}). It has been proved in \cite{BKT} that $K\in\R$ and, assuming that the additivity of the Lebesgue measure is $>\aleph_2$, $K$ is not an image of any Valdivia compact (see Remark 4.16 in \cite{BKT}).

Given a compact space $X$, we shall say that $p$ is a {\em cutting P-point} in $X$ if $X=A\cup B$, where $A,B$ are closed sets such that $A\cap B = \sn p$ and $p$ is a P-point both in $A$ and $B$. Recall that $p$ is a P-point in a space $Y$ if $p$ is not isolated in $Y$ and $p$ is not in the closure of any sequence of closed sets contained in $Y\setminus\sn p$. We prove that no space in class \R~contains cutting P-points. This immediately gives examples of spaces in $\RC\setminus\R$. For instance, the space $\omega_1+1+\omega_1^{-1}$, obtained from two copies of $\omega_1+1$ by identifying the two points of uncountable character, does not belong to $\R$. On the other hand, the space $\omega_1+1+\omega^{-1}$, obtained from the disjoint union of $\omega_1+1$ and $\omega+1$ by identifying the maximal elements, does belong to \R~(it can be obtained from $\omega_1+1$ by using a countable sequence of retractions).
Let us mention here that class \R~is not stable under open maps. An example is described in \cite{KU}, it is a compact connected Abelian group of weight $\aleph_1$. Every compact group is an epimorphic (and therefore open) image of a product of metric compact groups (i.e. of a Valdivia compact space).

The last section of this work is devoted to linearly ordered spaces. We use the result on cutting P-points for characterizing scattered linearly ordered compacta in class \R~and we show that all of them belong to \RC. Denote by $\ray+1$ the one-point compactification of the {\em long ray} $\ray=[0,1)\cdot \omega_1$ endowed with the lexicographic order. It turns out that the only non-metrizable connected linearly ordered spaces in class \R~are, up to order isomorphism, $\ray+1$, the inverse of $\ray+1$ and the {\em long interval} $\longint$, i.e. the space obtained by ``gluing'' $\ray+1$ together with its inverse in a suitable way. We also notice that all of these three spaces are Valdivia compact.

We finish with a study of Valdivia compact linearly ordered spaces. We observe that all of them have weight $\loe\aleph_1$ and we find a universal order preserving pre-image for this class. It is a $0$-dimensional dense-in-itself space $K$ whose all nonempty clopen subsets are order isomorphic to $K$. Moreover, every infinite interval contains both a strictly increasing and a strictly decreasing sequence of length $\omega_1$. The natural two-to-one order preserving quotient of $K$ gives a connected linearly ordered space in class \RC~which is nowhere separable.

\section{Preliminaries}

All spaces are assumed to be completely regular. A ``map" means a ``continuous map" unless otherwise indicated. A {\em retraction} or a {\em right-invertible map} is a map $\map fXY$ such that there exists $\map gYX$ with $fg=\id_Y$; $g$ is called a {\em right inverse} of $f$.
In this case $f$ is a quotient map and $g$ is an embedding. If $Y\subs X$ and $g$ is the inclusion map, then we say that $f$ is a retraction {\em into} $X$ or that $f$ is an {\em internal retraction} and $Y$ is a {\em retract} of $X$.

Given a linearly ordered set $\pair X\loe$, we shall consider its order topology, which is the one generated by open intervals of $X$. Recall that $\pair X\loe$ is compact if and only if every subset of $X$ has the least upper bound (the supremum of the empty set is the minimal element of $X$). We denote by $1$ the singleton $\sn0$ treated as a linearly ordered set. Given two disjoint linearly ordered sets $X,Y$, we denote by $X+Y$ the linearly ordered set whose universe is $X\cup Y$ and the order $\loe$ is defined by extending the union of the orders of $X$ and $Y$ and adding the relation $x<y$ for every $x\in X$ and $y\in Y$. When the sets $X,Y$ are not disjoint, we define $X+Y$ in a similar way, using an isomorphic copy of $Y$ which is disjoint from $X$. The sum defined above is a special case of the {\em lexicographic sum $\sum_{i\in Q}P_i$ of $\sett{P_i}{i\in Q}$ along} $Q$, where $P_i$'s and $Q$ are posets. The universe of $\sum_{i\in Q}P_i$ is $\bigcup_{i\in Q}(P_i\times\sn i)$ and the order is defined as follows:
$$\pair xi \loe \pair yj \iff i<j\text{ or }(i=j\text{ and }x\loe y\text{ in }P_i).$$
In case where $P_i=P$ for every $i$, the above lexicographic sum is denoted by $P\cdot Q$ and it is called the {\em lexicographic product} of $P$ by $Q$.
Given a linearly ordered set $X$ with the order $\loe$, we denote by $X^{-1}$ the set $X$ with the reversed ordering $\goe$. 
A linearly ordered set is {\em scattered} if it does not contain an isomorphic copy of the rationals $\Qyu$. In general, this is stronger than being topologically scattered in the order topology ($\Qyu\cdot \Nat$ is an example), however in the case of compact linearly ordered spaces both notions coincide.

Hausdorff's Theorem says that the class of all scattered linear orderings is the smallest class that contains ordinals and which is closed under reversing the order and under lexicographic sums along ordinals.

The class of all ordinals is denoted by $\Ord$. The first infinite and the first uncountable ordinal are denoted, as usual, by $\omega$ and $\omega_1$ respectively.


\subsection{Inverse systems}

Let $\S=\invsys Xp\Sigma st$ be an inverse system, i.e. $\Sigma$ is a directed partially ordered set, $\map{p^t_s}{X_t}{X_s}$ for every $s,t\in\Sigma$, $s\loe t$ and $p^s_s=\id_{X_s}$, $p^t_s p^r_t = p^r_s$ whenever $s\loe t\loe r$. Mappings $p^t_s$ are called {\em bonding maps}. 
We denote by $p_s$ the projection from $\liminv\S$ onto $X_s$. In the language of category theory, $\liminv\S$ is a pair consisting of a topological space $X$ and a family of {\em projections} $\setof{p_s}{s\in \Sig}$ with the property that for every topological space $Y$ and for every collection of maps $\setof{g_s}{s\in\Sig}$ such that $\map{g_s}Y{X_s}$ and $p^t_s g_t=g_s$ whenever $s\loe t$, there exists a unique map $\map gYX$ such that $p_sg=g_s$ holds for every $s\in\Sig$. In case where the projections are obviously defined, we shall denote by $\liminv\S$ the space itself.
A typical description of $\liminv\S$ is 
$$\liminv\S=\bigsetof{x\in\prod_{s\in\Sigma}X_s}{(\forall\;s,t\in\Sigma)\;s<t\implies p^t_s(x(t))=x(s)}$$
and the topology is inherited from the product. Given a cofinal set $T\subs \Sig$, the family $\setof{\inv{p_s}U}{s\in T,\;U\subs X_s\text{ is open}}$ is an open base for $\liminv\S$; thus $\liminv\S = \liminv(\S\rest T)$, where $\S\rest T:=\invsys XpTst$.
An inverse system of the form $\S=\invsys Xp\lam\xi\eta$, where $\lam$ is an ordinal with the natural order, is called an {\em inverse sequence}. The sequence $\S$ is {\em continuous}, if for every limit ordinal $\delta<\lam$ the space $X_\delta$ together with projections $\sett{p^\delta_\xi}{\xi<\delta}$ is homeomorphic to $\liminv(\S\rest\delta)$, where $\S\rest \delta=\invsys Xp\delta\xi\eta$.
We shall use the fact that every inverse sequence can be refined to a (cofinal) subsequence of a regular length.
In this work, we consider mostly continuous inverse sequences of compact spaces with surjections.

We say that $\S=\invsys Xr\Sigma st$ is an {\em inverse system of retractions} or a {\em retractive inverse system} if $\S$ is an inverse system in which each bonding map $r^t_s$ is a retraction.
We say that $\seqof{i^t_s}{s\loe t}$ is a {\em right inverse} of $\S$ if
$r^t_s i^t_s=\id_{X_s}$, $i^s_s=\id_{X_s}$ and $i^u_t i^t_s = i_s^u$,
whenever $s\loe t\loe u$.

The following simple properties of retractive sequences were proved in \cite{KM} and also in \cite{BKT}.

\begin{lm}[{\cite[Lemma 3.1]{KM}}]\label{lipa1} Assume $X=\liminv \S$, where $\S=\invsys Xr\Sigma s{s'}$ is an inverse system with a right inverse $\seqof{i^t_s}{s\loe t}$. Then there exist mappings $\map {i_s}{X_s}X$ such that
$r_s i_s = \id_X$ and $i_t i^t_s = i_s$, whenever $s\loe t$, $s,t\in\Sigma$.
\end{lm}

It is clear that maps $i_s$ in the above lemma are uniquely determined by the right inverse of $\S$. So, whenever $\S=\invsys Xr\Sigma st$ is an inverse system with a right inverse $\seqof{i^t_s}{s\loe t}$, we shall use mappings $\map{i_s}{X_s}X$ refering to Lemma \ref{lipa1} implicitly.

\begin{lm}[{\cite[Lemma 3.2]{KM}}]\label{lipa2} Let $\kappa$ be an ordinal, let $X=\liminv\S$, where $\S=\invsys Xr\kappa\al\beta$ is a continuous inverse sequence such that $r_\al^{\al+1}$ is a retraction for every $\al<\kappa$. Then $\S$ has a right inverse.
\end{lm}

\begin{lm}[{\cite[Lemma 3.3]{KM}}]\label{brzim} Let $\S=\invsys Xr\Sigma st$ be an inverse system with a right inverse $\seqof{i_s^t}{s\loe t}$ and define $R_s=i_s r_s$ for every $s\in \Sigma$. Then
\begin{enumerate}
	\item[$(a)$] $s\loe t\implies R_s R_t=R_s=R_t R_s$.
	\item[$(b)$] $x=\lim_{s\in\Sigma}R_s(x)$ for every $x\in \liminv\S$.
\end{enumerate}
\end{lm}

\begin{lm}[{\cite[Lemma 3.4]{KM}}]\label{antybrzim} Assume $\seqof{R_s}{s\in\Sig}$ is a family of retractions of a compact space $X$ into itself such that $\Sig$ is a directed partially ordered set and conditions (a), (b) of Lemma \ref{brzim} hold. Then $X=\liminv\invsys XR\Sig st$, where $X_s=\img{R_s}X$ and $R^t_s=R_s\rest X_t$.
\end{lm}

Given an ordinal $\delta$, by an {\em inverse sequence of internal retractions} of a space $X$ we mean a sequence $\setof{r_\al}{\al<\delta}$ of internal retractions of $X$ such that $r_\al r_\beta = r_\al = r_\beta r_\al$ for every $\al<\beta<\delta$ and $x=\lim_{\al<\delta}r_\al(x)$ for every $x\in X$. By the above lemmas, such a sequence describes uniquely a retractive inverse sequence with limit $X$. Continuity of this sequence is equivalent to saying that $r_\gamma(x)=\lim_{\al<\gamma}r_\al(x)$ for every limit ordinal $\gamma<\delta$ and for every $x\in X$. 

Let us note that every $0$-dimensional compact space is the limit of an inverse {\em system} of finite spaces and such a system always has a right inverse:

\begin{prop} Assume $\S=\invsys Xr\Sigma st$ is an inverse system of finite metric spaces such that each bonding map $r^t_s$ is a surjection. Then $\S$ has a right inverse.
\end{prop}

\begin{pf} Fix a well ordering $\rloe$ on $X=\liminv\S$. Each projection $\map{r_s}X{X_s}$ induces a partition of $X$ into clopen sets. Let $i_s(r_s(x))$ be defined as the $\rloe$-minimal element of $r_s^{-1}r_s(x)$. If $s\loe t$ then the partition induced by $r_t$ refines the one induced by $r_s$. More precisely, $r_t^{-1}r_t(x)\subs r_s^{-1}r_s(x)$, whenever $s\loe t$ and $x\in X$. Thus, setting $i^t_s(p)=r_t(i_s(x))$ we obtain a map $\map{i^t_s}{X_t}{X_s}$ which is a right inverse of $r^t_s$. Finally, if $s\loe t\loe u$ then $i_s^u=i_t^u i_s^t$.
\end{pf}

\subsection{Elementary substructures and quotients}

Given a regular cardinal $\chi$, we denote by $H(\chi)$ the class of all sets which are hereditarily of cardinality $<\chi$. We shall consider elementary substructures of the structure $\pair{H(\chi)}{\in}$, where $\chi$ is an uncountable (regular) cardinal. An important fact is that all countable elementary substructures of $\pair{H(\chi)}{\in}$ form a closed and cofinal set in $\pair{\dpower{H(\chi)}{\aleph_0}}\subs$. More specifically, the union of a chain of elementary substructures is an elementary substructure and, by the L\"owenheim-Skolem Theorem, every countable subset of $H(\chi)$ is contained in a countable elementary substructure of $H(\chi)$.
The fact that $M$ is an elementary substructure of $\pair{H(\chi)}{\in}$ will be denoted by $M\rloe\pair{H(\chi)}{\in}$ or just $M\rloe H(\chi)$. 
For sample applications of elementary substructures in topology we refer to \cite{Dow} and, in the context of Valdivia compacta, to \cite{KM}. 

Let $K$ be a compact space and let $C(K)$ denote its space of all real-valued continuous functions. Let $\chi>\aleph_0$ be such that $K\in H(\chi)$ and fix an elementary substructure $M$ of $\pair{H(\chi)}{\in}$ such that $K\in M$. By elementarity, $C(K)\cap M\nnempty$ and one can define an equivalence relation $\sim_M$ by
$$x\sim_M y\iff (\forall\;\phi\in C(K)\cap M)\;\phi(x)=\phi(y).$$
Let $K\by M$ denote the quotient $K\by{\sim_M}$ and let $q^M_K$ denote the quotient map. This construction has been first considered by Bandlow \cite{Ba91a, Ba91, Ba94} and used for studying some classes of compacta and related Banach spaces. We shall use a characterization of Corson compacta in terms of countable elementary substructures, proved by Bandlow in \cite{Ba91}.

The following statement is a consequence of \cite[Lemma 5]{Ba91} and \cite[Lemma 2.4]{KM}.

\begin{lm}\label{wnegfpjg} Let $X$ be Valdivia compact suitably embedded in $\cube\kappa$, i.e. $X=\cl(\Sig(\kappa)\cap X)$. Let $\chi>\kappa$ and let $M$ be an elementary substructure of $H(\chi)$ such that $X\in M$. Further, let $S=\kappa\cap M$. Then 
\begin{enumerate}
	\item[(a)] $\sim_M$ is the same as the equivalence relation induced by the map $\map{r_S}X{\cube\kappa}$ defined by $r_S(x)=x\cdot\chi_S$, where $\chi_S\in\cube \kappa$ denotes the characteristic function of $S$. 
	\item[(b)] $r_S$ is an internal retraction, i.e. $\img{r_S}X \subs X$ and $\img{r_S}X$ is Valdivia compact.
\end{enumerate}
\end{lm}


\section{Classes \RC~and \R}

As mentioned before, class $\RC$ is defined to be the smallest class that contains all metric compacta and which is closed under continuous images and inverse limits of transfinite sequences of retractions. Class $\RC$ can also be defined recursively as $\RC=\bigcup_{\xi\in\Ord}\RC^\xi$, where $\RC^0$ is the class of all metric compacta and $\RC^\beta$ consists of all compact spaces $X=\img f{\liminv\S}$, where $\S=\invsys Xr\kappa\al\beta$ is a continuous inverse sequence such that each $r^{\al+1}_\al$ is a retraction (see Lemma \ref{lipa2}) and each $X_\al$ belongs to $\bigcup_{\xi<\beta}\RC^\xi$. Given $X\in \RC$ denote by $\rk_{\RC}(X)$ the minimal $\beta$ such that $X\in\RC^\beta$. We call $\rk_\RC(X)$ the {\em \RC-rank of} $X$.

Class \R~has been defined in \cite{BKT} as the smallest class of spaces which contains all metric compacta and which is closed under limits of continuous inverse sequences of retractions. There is a natural hierarchy on $\R$. Denote by $\R^0$ the class of metric compacta and define $\R^\beta$ as the class of all spaces of the form $\liminv\S$, where $\S=\invsys Xr\kappa\xi\eta$ is a continuous inverse system such that each $r^{\xi+1}_\xi$ is a retraction and each $X_\xi$ belongs to $\bigcup_{\al<\beta}\R^\al$. Then $\R=\bigcup_{\xi\in\Ord}\R^\xi$. Given $X\in \R$ define $\rk_\R(X)$ as the minimal $\beta$ such that $X\in \R^\beta$. We call $\rk_\R(X)$ the {\em \R-rank of} $X$.

Clearly $\R\subs\RC$ and $\rk_\RC(X)\loe\rk_\R(X)$, whenever $X\in \R$.
Every Valdivia compact space belongs to class \R, since it can be decomposed into a retractive sequence of Valdivia compacta of smaller weight (see \cite{AMN}).
We shall prove more, namely that retracts of Valdivia compacta belong to \R~(see Theorem \ref{weiogeigf} below).

Both classes \RC~and \R~are closed under typical operations.
A covariant functor on topological spaces is {\em continuous} if it preserves limits of arbitrary inverse sequences.

\begin{prop}\label{funktor1} Assume $F$ is a continuous covariant functor on compact spaces.
\begin{enumerate}
  \item[(a)] If $F(X)\in\R$ for every compact metric space $X$ then $F$ preserves class \R. 
  \item[(b)] If $F(X)\in\RC$ for every compact metric space $X$ and $F(f)$ is a surjection whenever $f$ is a surjection, then $F$ preserves class $\RC$.
\end{enumerate}
\end{prop}

\begin{pf} (a) Define $\Kay=\setof{X\in\R}{F(X)\in\R}$. By assumption, $\Kay$ contains all metric compacta. By the continuity of $F$, $\Kay$ is closed under inverse limits of retractions (note that $F(r)$ is a retraction whenever $r$ is a retraction). Thus $\Kay=\R$.

(b) Define $\El=\setof{X\in\RC}{F(X)\in\RC}$. Then $\El$ contains all metric compacta and is closed both under limits of continuous retractive inverse sequences and under continuous images. The latter follows from the fact that $F$ preserves surjections. Thus $\El=\RC$.
\end{pf}

\begin{prop}\label{sgsdgsdgseewgs} Classes $\R$ and $\RC$ are stable under arbitrary products and one-point compactifications of disjoint topological sums. \end{prop}

\begin{pf} We give the proof for class \R; the case of class \RC~is the same, because  the functors considered below preserve surjections. 

An infinite product is the limit of an inverse sequence of smaller products, where the bonding mappings are projections.
Thus, it is enough to show that $X\times Y\in\R$ whenever $X,Y\in\R$.
Given a compact space $Y$, define $F_Y(X)=X\times Y$ and $F_Y(f)=f\times\id_Y$. Then $F_Y$ is a continuous covariant functor on compact spaces. Fix $Y\in\R$. If $Y$ is compact metric, then $F_Y(X)$ is compact metric for every compact metric space $X$.
Thus, by Proposition \ref{funktor1}, $F_Y(X)\in\R$ for every $X\in \R$. Now, if $Y$ is any space in \R, then $F_Y(X)=F_X(Y)\in\R$ whenever $X$ is compact metric, so again by Proposition \ref{funktor1}, $F_Y(X)\in\R$ for every $X\in \R$.

Now let $X$ be the one-point compactification of the disjoint sum $\bigoplus_{\xi<\kappa} X_\xi$, i.e. $X=\sn\infty\cup\bigcup_{\xi<\kappa} X_\xi$, where $X_\xi\cap X_\eta=\emptyset$ whenever $\xi\ne\eta$. Assume first that $\kappa$ is infinite. For each $\al<\kappa$ define $\map {r_\al}XX$ by $r_\al(x)=x$ if $x\in\bigcup_{\xi<\al}X_\xi$ and $r_\al(x)=\infty$ otherwise. It is straight to check that $\setof{r_\al}{\al<\kappa}$ is a continuous inverse sequence of internal retractions of $X$. By Lemma \ref{antybrzim}, $X=\liminv\invsys Yr\kappa\al\beta$, where $Y_\al=\img{r_\al}X$ is the one-point compactification of $\bigoplus_{\xi<\al}X_\xi$ and $r^\beta_\al=r_\al\rest Y_\beta$. By induction, it remains to prove that $X\oplus Y\in \R$ whenever $X,Y\in\R$. Given a compact space $Y$ define $G_Y(X)=X\oplus Y$ and $G(f)=f\oplus\id_Y$. Then $G_Y$ is a continuous covariant functor on compact spaces. The same argument as above shows that $G_Y(X)\in\R$ for every $X,Y\in \R$.
\end{pf}

It has been already mentioned that class \R~is not stable under open maps, see \cite{KU}. We do not know whether class \R~is stable under retractions. We also do not know whether every space in class \RC~is a continuous image of some space from class \R.
Let us finish this section with a model-theoretic type of stability. 

Let $X$ be a compact space in a ZFC model $M$ and let $\bbL=\closed^M(X)$, the lattice of closed subsets of $X$ defined in $M$. Then $X$ is naturally homeomorphic to the space of all ultrafilters over $\bbL$ (in fact all ultrafilters are principal).
Now let $N\sups M$ be another model of ZFC and let $X^N$ denote the space of all ultrafilters over $\bbL$, defined in $N$ ($\bbL$ is still a lattice of sets in $N$). In many cases, there are new ultrafilters, therefore usually $X^N\ne X$, although always $X\subs X^N$. We call $X^N$ the {\em interpretation of $X$ in $N$}. It can be shown easily that we obtain the same space, taking any sublattice of $\bbL$ which forms a closed base of $X$ in $M$. 
Given a continuous map of compact spaces $\map fXY$ in $M$, there is a unique map $\map{f^N}{X^N}{Y^N}$ which extends $f$, i.e. such that $\inv{(f^N)}a=\inv fa$ for every $a\in\closed^M(Y)$. Since $(fg)^N=f^N g^N$, this defines a functor from the class of compact spaces in $M$ into the class of compact spaces in $N$. 

Given a class of compact spaces $\Kay$, we say that $\Kay$ is {\em absolute} if for every two ZFC models $M\subs N$ and $X\in M$ such that $M\models$ ``$X\in\Kay$", we also have $N\models$ ``$X^N\in\Kay$", where $X^N$ denotes the interpretation of $X$ in $N$. To be formal, this defines {\em upward absoluteness}, however no class containing all metric compacta is downward absolute, unless it consists of all compact spaces (take a space $K$ not in the class and extend the universe by collapsing the weight of $K$ to $\aleph_0$; then the interpretation of $K$ in the extension is a second countable compact space). Absoluteness of some classes of compacta with respect to forcing extensions has been studied by Bandlow in \cite{Ba92}.

\begin{prop} Class \R~is absolute. More precisely: if $M\subs N$ are models of ZFC, $\al$ is an ordinal in $M$ and $X$ is a compact space in $M$ such that $M\models$ ``$\rk_\R(X)\loe\al$" then we also have $N\models$ ``$\rk_\R(X^N)\loe\al$". The same statement holds for class \RC.
\end{prop}

\begin{pf} Induction on $\al$ (ordinals in $M$). The statement is true if $\al=0$ since being second countable is absolute. Fix $X\in M$ such that $M\models$ ``$\rk_\R(X)=\al$" and let $X=\liminv\S$, where $\S=\invsys Xr\kappa\xi\eta\in M$ is a retractive continuous inverse sequence and $M\models$ ``$\rk_\R(X_\xi)<\al$" for every $\xi<\kappa$. 
By inductive hypothesis, $\rk_\R(X_\xi^N)<\al$ in $N$. Now observe that, setting $Y_\xi=X_\xi^N$ and $R_\xi^\eta=(r_\xi^\eta)^N$, we obtain a retractive inverse sequence $\S^N=\invsys YR\kappa\xi\eta$ in $N$. To see that $\S^N$ is continuous, notice that the functor $f\mapsto f^N$ preserves inverse limits. Indeed, any inverse system of quotient maps can be translated to a dual inductive system of lattices (the lattices of closed sets) and its inverse limit translates to the inductive limit. Now, the inductive limit is the same in any extension of the universe and the inverse limit is uniquely determined by the inductive limit of the dual system of lattices. It follows that $X^N=\liminv\S^N$ in $N$ and therefore $N\models$ ``$\rk_\R(X^N)\loe\al$".

It is clear that the above arguments can be adapted to show the absoluteness of class \RC. \end{pf}

\section{Main results}

In this section we prove the announced results on classes \R~and \RC.

\subsection{Retracts of Valdivia compacta}

\begin{tw}\label{weiogeigf} Retracts of Valdivia compact spaces belong to class \R.
\end{tw} 

\begin{pf} Induction on the weight of the retract. Assume $\kappa>\aleph_0$ is such that $X\in\R$ whenever $X$ is a retract of a Valdivia compact and $\w(X)<\kappa$. Fix a retraction $\map fYX$ such that $Y$ is Valdivia and $\w(X)=\kappa$. 

Fix a cardinal $\chi>\kappa$, so that $f\in H(\chi)$ and fix $N\rloe H(\chi)$ which contains a fixed base of $X$ and such that $f\in N$ and $|N|=\kappa$. Then $X\by N=X$.
Let $\sett{M_\al}{\al<\kappa}$ be a continuous chain of elementary substructures of $N$ with $\bigcup_{\al<\kappa}M_\al=N$, $f\in M_0$ and $|M_\al|<\kappa$ for every $\al<\kappa$. Let $q_\al$ denote the quotient map $X\to X\by M_\al$. Further, let $X_\al=\img{q_\al}X$ and for every $\al<\beta$ let $q^\beta_\al$ be the unique map such that $q_\al= q^\beta_\al q_\beta$ holds. Then $\S=\invsys Xq\kappa\al\beta$ is a continuous inverse sequence with limit $X$. It suffices to show that $\setof{X_\al}{\al<\kappa}\subs\R$ and that each $q_\al$ is a retraction.

Fix $\al<\kappa$. Let $r_\al$ denote the quotient map $Y\to Y\by M_\al$ and let $f_\al$ be the unique map such that the diagram
$$\xymatrix{
{Y} \ar[r]^-{f}\ar[d]_{r_\al} & {X}\ar[d]^{q_\al}\\
{Y\by M_\al} \ar[r]^-{f_\al} & {X_\al}
}
$$
commutes. By Lemma \ref{wnegfpjg}, $r_\al$ is a retraction and $Y\by M_\al$ is Valdivia compact. We claim that $f_\al$ is a retraction.

We may assume that $f$ is an internal retraction, i.e. $X\subs Y$ and $f\rest X=\id_X$. Let $Z=\img{r_\al}X\subs Y\by M_\al$. Given $x\in X$, we have $q_\al(x) = q_\al f(x) = f_\al r_\al(x)$. Observe that $r_\al$ and $q_\al$ induce the same equivalence relation $\sim_{M_\al}$ on $X$. This shows that $f_\al\rest Z$ is one-to-one. Clearly $\img{f_\al}Z=q_\al\img fY=X_\al$. 
Thus $f_\al$ is right-invertible. By the induction hypothesis, $X_\al=\img{f_\al}{Y\by M_\al}\in\R$. Finally, the composition $q_\al f = f_\al r_\al$ is a retraction, hence so is $q_\al$.
\end{pf}

\subsection{The dichotomy}

We start with an auxiliary result on retractive inverse sequences. Recall that given an inverse sequence of internal retractions, we may represent its limit as the closure of the union of all the spaces from the sequence, where the projections are internal retractions (see Lemma \ref{brzim}).

\begin{lm}\label{faoiwqurjfoar}
Let $\S=\invsys Xr\omega nm$ be an inverse sequence of internal retractions of compact spaces with $X=\liminv\S$, represented in such a way that all projections $\map{r_n}X{X_n}$ are internal retractions.
Let $\map fXY$ be a surjection and let,
for each $\ntr$, $\Bee_n$ be an open base for $\img f{X_n}$ which is closed under finite unions. Then the family 
$$\Cee=\setof{Y\setminus \img f{X\setminus \inv{(fr_n)}B}}{\ntr,\;B\in\Bee_n}$$
is an open base for $Y$.
In particular, $\w(Y) = \sup_{\ntr}\w(\img f{X_n})$.
\end{lm}

\begin{pf}
Fix $y\in Y$ and its neighborhood $U$. Find $m\in\nat$ and an open set $W\subs X_m$ such that
\begin{equation}
f^{-1}(y)\subs\inv {r_m}W\subs \inv fU.
\tag{1}\end{equation}
Let $F=\img f{X\setminus \inv{r_m}W}$. Then $y\notin F$, therefore $U_1=Y\setminus F$ is a neighborhood of $y$, contained in $U$. 
Now find $n\goe m$ and an open set $W_1\subs X_n$ such that
$$f^{-1}(y)\subs \inv{r_n}{W_1}\subs \inv f{U_1}.$$
Thus $\img{r_n}{f^{-1}(y)}\subs W_1$ and $W_1\subs \inv f{U_1}$. The latter inclusion follows from the fact that $W_1\subs \inv{r_n}{W_1}$. Thus $\img {fr_n}{f^{-1}(y)}\subs U_1$.
Using compactness and the fact that $\Bee_n$ is closed under finite unions, we can find $B\in\Bee_n$ such that
$$\img {fr_n}{f^{-1}(y)}\subs B\subs U_1.$$
We have 
$$\img{r_n}{f^{-1}(y)}\subs \inv fB\subs \inv{r_m}W,$$
where the latter inclusion follows from the equality $U_1=Y\setminus \img f{X\setminus\inv{r_m}W}$.
Thus
\begin{equation}
f^{-1}(y) \subs \inv{r_n}{\inv f B}\subs \inv{r_n}{\inv{r_m}W}\subs \inv{(r_mr_n)}W=\inv{r_m}W
\tag{2}\end{equation}
Hence, using (1) and (2), we obtain $f^{-1}(y)\subs\inv{r_n}{\inv fB}\subs \inv fU$.
Finally, set $V=Y\setminus \img f{X\setminus\inv{(fr_n)}B}$. Then $V\in\Cee$ and $y\in V\subs U$.
\end{pf}

Let us recall a model-theoretic characterization of Corson compacta, due to Bandlow \cite{Ba91}. A compact space $X$ is Corson compact if and only if for a big enough cardinal $\chi$ there exists a closed and unbounded family $\Em\subs\dpower{H(\chi)}{\aleph_0}$ of elementary substructures of $H(\chi)$ such that for every $M\in\Em$ the quotient map $\map{q_M}X{X\by M}$ is one-to-one on $\cl(X\cap M)$. Here, {\em closed and unbounded} means: closed under the unions of countable chains and cofinal in $\pair{\dpower{H(\chi)}{\aleph_0}}{\subs}$. In fact, if $X$ is Corson, then the above property holds for {\em every} countable elementary substructure of $H(\chi)$ which ``knows" an embedding of $X$ into a $\Sig$-product. See \cite{Ba91} for more details.

\begin{tw}\label{wetijgfjsfjipads} Let $X\in\RC$. Then either $X$ is Corson compact or else $X$ contains a copy of the linearly ordered space $\omega_1+1$.
\end{tw}

\begin{pf} Every second countable compact space is Corson, therefore the above dichotomy holds for spaces of \RC-rank $0$. Fix an ordinal $\beta>0$ and assume that the above statement is true for spaces of \RC-rank $<\beta$. Fix $X\in\RC$ such that $\rk_\RC(X)=\beta$. Let $\S=\invsys Yr\kappa\xi\eta$ be a continuous inverse system of internal retractions such that $X=\img fY$ for some map $\map fYX$, where $Y=\liminv\S$ and $\rk_\RC(Y_\xi)<\beta$ for every $\xi<\kappa$. We assume that $\kappa$ is a regular cardinal. Let $X_\xi=\img f{Y_\xi}$ for $\xi<\kappa$ and let $X_\kappa=X$. If $X_\xi$ contains a copy of $\omega_1+1$ for some $\xi<\kappa$, then there is nothing to prove, so assume that $X_\xi$ is Corson compact for each $\xi<\kappa$.
Suppose first that
\begin{enumerate}
	\item[(*)] $\kappa>\aleph_0$ and $X\ne\bigcup_{\xi<\kappa}X_\xi$.
\end{enumerate}
Then $Y\ne\bigcup_{\xi<\kappa}Y_\xi$. Fix $y\in Y$ so that $f(y)\notin\bigcup_{\xi<\kappa}X_\xi$.
Recall that $y=\lim_{\xi\to\kappa}r_\xi(y)$. Also $r_\delta(y)=\lim_{\xi<\delta}r_\xi(y)$, by the continuity of the sequence $\S$. Thus the map $\map\phi{\kappa+1}Y$, defined by $\phi(\xi)=fr_\xi(y)$ and $\phi(\kappa)=f(y)$, is continuous. By assumption, $\phi(\xi)\ne f(y)$ for every $\xi<\kappa$, therefore the sequence $\sett{\phi(\xi)}{\xi<\kappa}$ does not stabilize. Find a closed cofinal set $C\subs\kappa$ such that $\phi\rest C$ is one-to-one. Then $Z=\img \phi{C\cup\sn\kappa}\subs X$ is homeomorphic to the linearly ordered space $\kappa+1$. Since $\kappa>\aleph_0$, this shows that $X$ contains $\omega_1+1$.

Now suppose that (*) does not hold. We claim that $X$ is Corson. For this aim we use Bandlow's characterization. Fix a big enough cardinal $\chi$ and a countable $M\rloe H(\chi)$ such that $f,\S\in M$.

Let $\delta=\sup(\kappa\cap M)$ (if $\kappa=\aleph_0$ then of course $\delta=\omega\subs M$ and $X_\delta=X$). Define $X_\xi^M=\cl(X_\xi\cap M)$ and $Y_\xi^M=\cl(Y_\xi\cap M)$ for $\xi<\kappa$. Let $\map {q_M}X{X\by M}$ denote the quotient map induced by $M$. In order to show that $q_M$ is one-to-one on $X^M:=\cl(X\cap M)$, it suffices to find a base for $X^M$ which is contained in $M$.

Notice that $X\cap M\subs X_\delta$. Indeed, if $\kappa=\aleph_0$ then $X_\delta=X$; otherwise $X=\bigcup_{\xi<\kappa}X_\xi$ and therefore by elementarity $X\cap M\subs\bigcup_{\xi\in\kappa\cap M}X_\xi\subs X_\delta$.

Now observe that $\img f{Y_\xi^M}=X_\xi^M$ for every $\xi\in \kappa\cap M$. 
Indeed, if $x\in X_\xi\cap M$ and $\xi \in M$ then by elementarity $x=f(y)$ for some $y\in Y_\xi\cap M$, which shows that $X_\xi\cap M\subs \img f{Y_\xi\cap M}$. Clearly $\img f{Y_\xi\cap M}\subs X_\xi\cap M$.

Note that $\img f{Y^M}=X^M$, because $f\in M$. On the other hand, $Y^M=Y_\delta^M=\cl(\bigcup_{\xi\in\kappa\cap M}Y_\xi^M)$. Let $\setof{\xi_n}{\ntr}\subs \delta\cap M$ be an increasing sequence such that $\delta=\sup_{\ntr}\xi_n$. Let $\S_0=\invsys Zp\omega nm$, where $Z_n=Y_{\xi_n}^M$, $p_n=r_{\xi_n}\rest Z_n$. We apply Lemma \ref{faoiwqurjfoar} to find a suitable base for $X^M$. 
Fix $\ntr$. Then $X_{\xi_n}\in M$ and by assumption it is Corson compact. Thus $M$ induces a quotient map which is one-to-one on $X_{\xi_n}^M$. In particular, the family $$\Bee_n=\setof{U\in M}{U\text{ is an open subset of }X_{\xi_n}}$$
forms a base for $X_{\xi_n}^M$ (more precisely: $\setof{U\cap X_{\xi_n}^M}{U\in\Bee_n}$ is an open base for $X_{\xi_n}^M$). Clearly, $\Bee_n$ is closed under finite unions. Lemma \ref{faoiwqurjfoar} says that
$$\Cee=\setof{X^M\setminus \img f{Y^M\setminus \inv{(fp_n)}B}}{\ntr,\;B\in\Bee_n}$$
is an open base for $X^M=\img f{Y^M}$.
Let $\Cee'=\setof{X\setminus \img f{Y\setminus \inv{(fp_n)}B}}{\ntr,\;B\in\Bee_n}$. Then $\Cee'\subs M$ and $\setof{C\cap X^M}{C\in\Cee'}=\Cee$. Thus we have shown that $M$ contains a family which forms an open base for $X^M$. Since $M$ was arbitrary, this implies that $X$ is Corson compact.
\end{pf}

\begin{wn}\label{owijep} Countably tight spaces in class \RC~are Corson compact. \end{wn}

Theorem \ref{wetijgfjsfjipads} gives immediately simple examples of compact spaces which do not belong to class \RC. For example, Mr\'owka's space, which can be defined as the Stone space of the Boolean algebra generated by $\fin\omega\cup\Aaa\subs\Pee(\omega)$, where $\Aaa$ is an uncountable almost disjoint family, is a scattered space of height $2$ and it is countably tight, separable and not metrizable. Thus, it is not Corson and therefore not a member of class \RC.

\subsection{Cutting P-points}

Recall that a point $p\in X$ is called a {\em P-point} if it is not isolated in $X$ and for every sequence $\ciag u$ of neighborhoods of $p$ it holds that $p\in\Int(\bigcap_{\ntr}u_n)$. We shall say that $p\in X$ is a {\em cutting P-point} if there exist closed subspaces $A,B$ of $X$ such that $A\cup B=X$, $A\cap B=\sn p$ and $p$ is a P-point both in $A$ and $B$.

A typical example of a cutting P-point is the complete accumulation point of the linearly ordered space $\omega_1+1+\omega_1^{-1}$. 

\begin{tw}\label{ewjgpjgpjeg} 
No space in class \R~contains cutting P-points.
\end{tw}

\begin{pf}
Suppose the above statement is not a theorem and fix a model of ZFC which contains a counterexample $K$, witnessed by closed sets $A,B\subs K$ with $A\cap B=\sn p$. Let $\kappa=\w(K)^+$ and extend the universe by using the natural forcing which collapses $\kappa$ to $\aleph_1$. Then the interpretation of $K$ in the extension is still a counterexample, because it belongs to $\R$ and there are no ``new" countable sequences, so $p$ remains a P-point both in $A$ and $B$. Thus, working in a suitable model of ZFC, we may assume that $\w(K)=\aleph_1$. We may also assume that $\rk_\R(K)$ is minimal possible.

Let $\setof{u_\al}{\al<\omega_1}$ denote a basis of $p$ in $K$ such that $\al<\beta\implies \cl u_\beta\subs u_\al$. Such a sequence exists, because $p$ is a P-point in $K$ and $\chi(p,K)=\aleph_1$. The assumptions on $A,B$ mean in particular that both sets $u_\al\setminus B$, $u_\al\setminus A$ are nonempty for every $\al<\omega_1$.

Now suppose $\setof{r_\al}{\al<\kappa}$ is a sequence of internal retractions in $K$, such that $K_\al=\img{r_\al}K$ has $\R$-rank $<\rk_\R(K)$ for every $\al<\kappa$. Assuming $\kappa$ is regular, we deduce that $\kappa\in\dn{\aleph_0}{\aleph_1}$.
Recall that $\bigcup_{\al<\kappa}K_\al$ is dense in $K$.

Suppose $\kappa=\aleph_0$. For each $\al<\omega_1$ find $n(\al)<\omega$ with $u_\al\cap A\cap K_{n(\al)}\nnempty\ne u_\al\cap B\cap K_{n(\al)}$. This is possible, because $p$ is an accumulation point of both $A$ and $B$ and therefore $u_\al\setminus B$ and $u_\al\setminus A$ are nonempty open sets. Now there is $\ntr$ such that $\setof{\al<\omega_1}{n(\al)=n}$ is uncountable. Then $p$ is an accumulation point of both $A\cap K_n$ and $B\cap K_n$ which implies that $K_n$ is a counterexample to the theorem. This contradicts the minimality of $\rk_\R(K)$. It follows that $\kappa=\aleph_1$. 

The same argument as above, using the minimality of $\rk_\R(K)$, shows that either
\begin{equation}
(\forall\;\al<\omega_1)(\exists\;\beta<\omega_1)\; K_\al\cap u_\beta\setminus B =\emptyset
\tag{1}\end{equation}
or the same holds for $A$ in place of $B$. 
Clearly, we may assume that either
\begin{equation}
(\forall\;\al<\omega_1)\; r_\al(p)\in B
\tag{2}\end{equation}
or the same holds for $A$ in place of $B$.
Thus, without loss of generality, we may assume that both (1) and (2) hold (interchanging the roles of $A$ and $B$, if necessary). 
Indeed, if $p\notin\bigcup_{\al<\omega_1}K_\al$ then (1) holds both for $A$ and $B$ and if $p\in K_\al$ for some $\al<\omega_1$ then $r_\beta(p)=p\in A\cap B$ for $\beta\goe \al$, therefore (2) holds for both $A$ and $B$.

Define $\ovr u_\al = \bigcap_{\xi<\al}\cl u_{\xi+1}$. Then $u_\al\subs \ovr u_\al$ and $\bigcap_{\al<\omega_1}\ovr u_\al=\sn p$, thus every neighborhood of $p$ contains some $\ovr u_\al$. Moreover, for a limit ordinal $\delta<\omega_1$ we have $\ovr u_\delta=\bigcap_{\xi<\delta}u_\xi=\bigcap_{\xi<\delta}\cl u_\xi$.

Fix $\al<\omega_1$. Observe that by (1) and (2) the set $B\cap K_\al$ is a neighborhood of $r_\al(p)$ in $K_\al$. By the continuity of $r_\al$, there is $\xi(\al)<\omega_1$ such that 
\begin{equation}
\img{r_\al}{\ovr u_{\xi(\al)}}\subs B.
\tag{3}\end{equation}
On the other hand, the set $u_\al\setminus B$ is open and nonempty, therefore there is $\eta(\al)>\al$ such that
\begin{equation}
K_{\eta(\al)}\cap u_\al\setminus B\nnempty.
\tag{4}\end{equation}
Fix a limit ordinal $\delta<\omega_1$ such that $\xi(\al)<\delta$ and $\eta(\al)<\delta$ whenever $\al<\delta$. Then, by (3), we have $\img{r_\al}{\ovr u_\delta}\subs B$ for every $\al<\delta$ and hence $\img{r_\delta}{\ovr u_\delta}\subs B$, because $r_\delta(x)=\lim_{\al<\delta}r_\al(x)$ for every $x\in K$. 

Now for each $\al<\delta$ fix $q_\al\in K_\delta\cap u_\al\setminus B$, which exists by (4). Let $q$ be an accumulation point of the sequence $\setof{q_\al}{\al<\delta}$. Then $q\in K_\delta\cap A$ and $q\in \bigcap_{\al<\delta}u_\al=\ovr u_\delta$. Furthermore, $q\ne p$, because $p$ is a P-point in $A$. Hence $q\notin B$ and $r_\delta(q)=q$, which contradicts the fact that $\img{r_\delta}{\ovr u_\delta}\subs B$.
\end{pf}

An immediate consequence of the above result is that $\omega_1+1+\omega_1^{-1}\notin \R$, where $\omega_1^{-1}$ denotes $\omega_1$ with the reversed ordering.

\begin{wn} Assume $K\in \R$ is $0$-dimensional and contains at least two P-points. Then there is a two-to-one map $\map fKL$ onto $L$ such that $L\notin\R$.
\end{wn}

\begin{pf} Let $p,q$ be two distinct P-points in $K$. Let $A\subs K$ be clopen such that $p\in A$ and $q\in K\setminus A$. Let $L$ be the quotient of $K$ induced by identifying points $p,q$. Then $L$ contains a cutting P-point (the image of $p$). Thus $L\notin\R$, by Theorem \ref{ewjgpjgpjeg}.
\end{pf}

\section{Linearly ordered spaces}

In this section we discuss linearly ordered spaces in classes \RC~and \R~and we prove the results announced in the introduction. 
Recall that every linearly ordered Corson compact is metrizable, see Theorem IV.10.1 on page 181 in \cite{A}. Thus, Theorem \ref{wetijgfjsfjipads} implies that spaces like a double arrow or a compact Aronszajn line do not belong to class \RC. 


\subsection{Connected spaces}

Recall that the {\em long ray} is the linearly ordered space $\ray=\omega_1\cdot[0,1)$ (where $P\cdot Q$ denotes the lexicographic product of $P$ by $Q$). Note that $\ray+1$ is a connected linearly ordered compact space with exactly one point of uncountable character (namely, the maximal point). Moreover, each proper interval which does not contain the maximal point is separable. Denote by $\longint$ the unique linearly ordered space $K=[a,b]$ such that for every internal point $p\in (a,b)$ the interval $[p,b]$ is order isomorphic to $\ray+1$ and the interval $[a,p]$ is order isomorphic to the inverse of $\ray+1$.

\begin{prop}
The spaces $\ray+1$ and $\longint$ are Valdivia compact.
\end{prop}

\begin{pf}
Let $K_\al=\delta\cdot[0,1)+1$ and define $\map{r_\al}{\ray+1}{K_\al}$ by setting $r_\al(x)=x$ if $x\in K_\al$ and $r_\al(x)=\max(K_\al)$ otherwise. Then $\setof{r_\al}{\al<\omega_1}$ is a continuous inverse sequence of internal retractions in $\ray+1$ and each $K_\al$ is second countable. By Theorem \ref{KM}, this shows that $\ray+1$ is Valdivia compact. The case of $\longint$ is similar. 
\end{pf}

\begin{tw}\label{wejqwfjqpw} Let $K\in\R$ be a connected linearly ordered space. Then $K$ is order isomorphic to one of the following spaces: $1$, $[0,1]$, $\ray+1$, $(\ray+1)^{-1}$, $\longint$. In particular, $K$ is Valdivia compact.
\end{tw}

\begin{pf}
Suppose the statement is false and let $K$ be a counterexample of a minimal possible \R-rank. Fix an inverse sequence of internal retractions $\setof{r_\al}{\al<\kappa}$ in $K$ such that $\rk_\R(K_\al)<\rk_\R(K)$ for every $\al<\kappa$, where $K_\al=\img{r_\al}K$. We assume that $\kappa$ is a regular cardinal.
Since $K$ is connected, each $K_\al$ is an interval of the form $[a_\al,b_\al]$. Further, the sequence $\sett{a_\al}{\al<\kappa}$ is non-increasing and the sequence $\sett{b_\al}{\al<\kappa}$ is non-decreasing. Taking a cofinal set of indices and reversing the order of $K$ if necessary, we may assume that $\sett{b_\al}{\al<\kappa}$ is strictly increasing. Pick $p\in(a_0,b_0)$. By the minimality of $\rk_\R(K)$, each of the intervals $[a_\al,b_\al]$ is order isomorphic to one of the spaces listed in the theorem. It follows that $[p,b_\al]\subs(a_{\al+1},b_{\al+1})$ is order isomorphic to $[0,1]$. Hence, $\kappa\in\dn{\aleph_0}{\aleph_1}$. Let $K=[a,b]$. Then $[p,b]$ is order isomorphic either to $[0,1]$ (if $\kappa=\aleph_0$) or to $\ray+1$ (if $\kappa=\aleph_1$).
By the same argument, $[a,p]$ is order isomorphic either to $[0,1]$ or to the reversed $\ray+1$. This shows that $K$ is not a counterexample, a contradiction.
\end{pf}

\subsection{Scattered spaces}

We already know that the space $\omega_1+1+\omega_1^{-1}$ is not in \R, while $\omega_1+1+\omega^{-1}\in\R$. Below we characterize scattered linearly ordered spaces which are in \R. Let $K$ be a linearly ordered compact space. It is not hard to check that $p\in K$ is a cutting P-point in $K$ if and only if $p$ is a P-point both in $(\leftarrow,p]$ and $[p,\rightarrow)$. 

\begin{tw} Let $K$ be a scattered linearly ordered compact space. Then $K\in\R$ if and only if $K$ has no cutting P-points.
\end{tw}

\begin{pf} The necessity follows from Theorem \ref{ewjgpjgpjeg}. Fix a compact scattered linearly ordered space $K$ with no cutting P-points. 
We use the idea from the proof of Hausdorff theorem on the structure of scattered linearly ordered sets.
Define the following relation on $K$.
$$x\sim y\iff (\forall\;a,b\in[x,y])\; [a,b]\in\R.$$
We check that $\sim$ is an equivalence relation. Only transitivity requires an argument. So assume $x<y<z$ and $x\sim y$, $y\sim z$. Fix $a\in [x,y)$, $b\in(y,z]$. If $y$ is isolated from, say, right-hand side, then $(y,b]=[y',b]\in\R$, $[a,y]\in\R$ and hence $[a,b]=[a,y]\oplus [y',b]\in\R$, by Proposition \ref{sgsdgsdgseewgs}. 
Otherwise, reversing the order if necessary, we may assume that $y=\inf\ciag y$, where $\ciag y$ is strictly decreasing. Assume $y_0=b$. Since $K$ is compact and $0$-dimensional, we may also assume that each $y_n$ has an immediate successor $z_{n-1}$, i.e. $y_n<z_{n-1}$ and $(y_n,z_{n-1})=\emptyset$. Now
$$[y,b]=\sn y\cup\bigcup_{\ntr}[z_n,y_n]$$
and the intervals $[z_n,y_n]$ are pairwise disjoint. Further, by assumption $[a,y]\in\R$ and hence, setting $K_n=[a,y]\cup\bigcup_{i<n}[z_i,y_i]$, we have that $K_n\in\R$, by Proposition \ref{sgsdgsdgseewgs}. Define $\map{r_n}{[a,b]}{K_n}$ by $r_n\rest K_n=\id_{K_n}$ and $r_n(x)=y$ for $x\in [a,b]\setminus K_n$. Then $\setof{r_n}{\ntr}$ is an inverse sequence of internal retractions in $[a,b]$. This shows that $[a,b]\in\R$ and completes the proof of the transitivity of $\sim$.

By definition, the equivalence classes of $\sim$ are convex. We show that they are closed. Fix $C\in K\by\sim$ and let $b=\sup C$. Fix $a\in C$. Let $\setof{b_\al}{\al<\kappa}\subs[a,b)$ be strictly increasing, continuous and such that $b=\sup_{\al<\kappa}b_\al$. Then $[a,b_\al]\in\R$, because $a,b_\al\in C$. Define $\map{r_\al}{[a,b]}{[a,b_\al]}$ by $r_\al(x)=\min\{x,b_\al\}$. Then $\setof{r_\al}{\al<\kappa}$ is a continuous inverse sequence of internal retractions of $[a,b]$. Hence $[a,b]\in\R$. This shows that $x\sim b$ for every $x\in C$, i.e. $b\in C$. By the same argument, $\inf C\in C$, i.e. $C$ is a closed interval.

Finally, it remains to show that there is only one $\sim$-equivalence class. Suppose not, and fix $C,D\in K\by\sim$, $\max C=c<d=\min D$. Then $c\not\sim d$ and in particular $(c,d)\nnempty$. Thus, there exists $E\in K\by\sim$ such that $C<E<D$, where $A<B$ means $\sup A<\inf B$. It follows that the quotient ordering of $K\by\sim$ is dense, i.e. $K\by \sim$ is dense-in-itself with the order topology. Note that the order topology of $K\by\sim$, induced by the quotient ordering, coincides with the quotient topology. Thus $K\by\sim$ is a dense-in-itself continuous image of $K$, which contradicts the fact that $K$ is scattered.
\end{pf}

\begin{wn} Let $K$ be a scattered compact linearly ordered space. Then there exists a scattered compact linearly ordered space $L\in\R$ which has a two-to-one continuous order preserving map onto $K$. In particular $K\in\RC$.
\end{wn}

\begin{pf}
Let $L=K\cdot2$, the lexicographic product of $K$ by $\{0,1\}$. Then $L$ is scattered, compact and it has no cutting P-points (every point of $L$ is isolated from at least one side). Thus $L\in \R$ and of course the natural projection $\map fLK$ is two-to-one and continuous.
\end{pf}

\subsection{Valdivia compacta}

\begin{prop}\label{bnjerg} Valdivia compact linearly ordered spaces have weight at most $\aleph_1$.
\end{prop}

\begin{pf} Suppose the above statement is not a ZFC theorem and fix a countable transitive model $M$ of ZFC which contains a counterexample $K$. Let $\poset\in M$ be the standard forcing notion which collapses $\aleph_1$ to $\aleph_0$. Let $G$ be a $\poset$-generic filter over $M$ and let $K^G$ denote the interpretation of $K$ in $M[G]$. Then $M[G]\models$ ``$\w(K^G)\goe\aleph_1$", because $\w(K)\goe\aleph_2$ in $M$. In order to get a contradiction, it suffices to show that $M[G]\models$ ``$\chi(K^G)\loe\aleph_0$", because then $K^G$ is a non-metrizable Corson compact linearly ordered space in $M[G]$ (being linearly ordered is absolute). 

We suppose that $M[G]\models$ ``$\chi(K^G)>\aleph_0$". Since $K^G$ is linearly ordered and $K$ is a dense subset of $K^G$, there is in $M[G]$ a strictly increasing or a strictly decreasing function $\map f{\aleph_2^M}K$ (recall that $\aleph_2^M=\aleph_1^{M[G]}$). Since $M\models$ ``$|\poset|=\aleph_1$", there is $S\in(\dpower{\aleph_2}{\aleph_2})^M$ such that $f\rest S\in M$. It follows that
$M\models$ ``$K$ contains a monotone sequence of length $\aleph_2$" and consequently $\chi(K)>\aleph_1$ in $M$.

We now work in $M$, showing that $\chi(K)\loe\aleph_1$, which gives a contradiction. The argument is very similar to the one which shows that $\omega_2+1$ is not Valdivia.

Suppose $\setof{p_\al}{\al<\omega_2}\subs K$ is strictly increasing and continuous, $p=\sup_{\al<\omega_2}p_\al$. Assume $K\subs[0,1]^\kappa$ so that $D=\Sigma(\kappa)\cap K$ is dense in $K$. For each $\al<\omega_2$ 
choose $q_\al\in (p_\al,p_{\al+2})\cap D$. Then $p=\lim_{\al<\omega_2}q_\al$ and $p$ is not in the closure of any countable subset of $\setof{q_\al}{\al<\omega_2}$. Thus $p\notin D$, because $D$ has a countable tightness. Let $S=\suppt(p)$. Then $|S|>\aleph_0$. Since $p=\lim_{\al<\omega_2}$, for each $s\in S$ there is $\al(s)<\omega_2$ such that $p_\al(s)>0$ for $\al\goe\al(s)$. Fix $T\in\dpower S{\aleph_1}$ and let $\beta=\sup_{s\in T}\al(s)$. Then $p_{\beta+1}(s)>0$ for every $s\in T$, i.e. $p_{\beta+1}\notin D$. This is a contradiction.
\end{pf}

We say that a space is {\em nowhere first countable} if its first countable subsets have empty interior. Note that given a linearly ordered space $X\in\RC$, the property of being nowhere first countable is equivalent to being nowhere separable. Indeed, if $J\subs X$ is a closed interval then $J\in\RC$; thus, if $J$ is first countable then, by Corollary \ref{owijep}, it is a linearly ordered Corson compact, hence second countable. On the other hand, if $p\in J$ does not have a countable base in $J$ then there exists a monotone sequence of length $\omega_1$ in $J$; the closure of such a sequence is, by compactness, homeomorphic to $\omega_1+1$.

The following two statements concerning linearly ordered compacta are well known. We give the proofs for the sake of completeness. Recall that, for compact linearly ordered spaces, every order preserving epimorphism is continuous. In particular, every order isomorphism is a homeomorphism.

\begin{prop} Every metrizable 0-dimensional dense-in-itself linearly ordered compact space is order isomorphic to the Cantor set $2^\omega$, endowed with the lexicographic ordering $<_{\lex}$ defined by $x<_{\lex}y$ iff $x(n)<y(n)$, where $n=\min\setof{k}{x(k)\ne y(k)}$.
\end{prop}

\begin{pf} It is clear that the Cantor set satisfies the above assumptions. Let $\pair X<$ and $\pair Y<$ be two spaces satisfying the above assumptions. Let $Q(X)$ denote the set of all points $x\in X$ such that $[x',x]=\dn {x'}x$ for some $x'<x$. Then $Q(X)$ is countable, because $X$ is metrizable. 
Observe that $\pair {Q(X)}<$ is a dense ordering, since $X$ is dense-in-itself and 0-dimensional. The same is true for $\pair{Q(Y)}<$, therefore both $Q(X)$, $Q(Y)$ are order isomorphic to the rationals (by Cantor's theorem).
Finally, any order isomorphism of $Q(X)$ onto $Q(Y)$ can be, by compactness, uniquely extended to an order isomorphism of $X$ and $Y$.
\end{pf}

\begin{prop} Let $\pair X<$ be a compact linearly ordered space of weight $\aleph_1$. Then $X$ is the limit of a continuous inverse sequence of metric linearly ordered compacta, where all the bonding maps, as well as all the projections, are order preserving surjections.
\end{prop}

\begin{pf} Fix a family $\Yu$ of open intervals of $X$ which forms a base for $X$ and so that $|\Yu|=\aleph_1$. Given $U,V\in \Yu$ with disjoint closures, choose a continuous order preserving function $f=f_{U,V}$ which separates $U$ from $V$. 
The existence of such a function $f$ is a consequence of the Urysohn Lemma. Indeed,
given a continuous function $\map h X[0,1]$ such that $U\subs h^{-1}(0)$, $V\subs h^{-1}(1)$ and assuming $\sup U=a<b=\inf V$, the conditions
$$f\rest (\leftarrow,a]=0,\quad f\rest [b,\rightarrow)=1\quad\text{and}\quad f(x)=\max h [a,x]\quad\text{for }a<x<b,$$
define a continuous order preserving function $f$ satisfying $f(\sup U)=0<1=f(\inf V)$.

Let $\Ef$ consist of all the functions $f_{U,V}$, where $U,V\in\Yu$ have disjoint closures. Then $|\Ef|=\aleph_1$ and for every $x<y$ in $X$ there exists $f\in \Ef$ with $f(x)<f(y)$. Let $\Ef=\bigcup_{\al<\omega_1}\Ef_\al$, where $\sett{\Ef_\al}{\al<\omega_1}$ is a continuous chain of sets such that $|\Ef_\al|=\aleph_0$ for every $\al<\omega_1$. Let $\sim_\al$ be the following relation:
$$x\sim_\al y \iff (\forall\;f\in\Ef_\al)\;f(x)=f(y).$$
It is clear that $\sim_\al$ is an equivalence relation and all equivalence classes of $\sim_\al$ are closed. Since $\Ef_\al$ consists of order preserving functions, the $\sim_\al$-equivalence classes are convex.

Let $X_\al=X\by{\sim_{\al}}$ and let $\map {q_\al}X{X_\al}$ be the quotient map. Then $X=\liminv \S$, where $\S=\invsys Xq{\omega_1}\al\beta$ and $q^\beta_\al$ is the unique map satisfying $q_\al = q^\beta_\al q_\beta$. Since $\sett{\Ef_\al}{\al<\omega_1}$ is continuous (i.e. $\Ef_\delta=\bigcup_{\xi<\delta}\Ef_\xi$ for every limit ordinal $\delta<\omega_1$), we infer that the inverse sequence $\S$ is continuous.

Finally, for every $\al<\omega_1$, the relation $\sim_\al$ induces a linear order $<_\al$ on $X_\al$ defined by
$$q_\al(x) <_\al q_\al(y) \iff x<y.$$
This is well defined, because the $\sim_\al$-equivalence classes are convex. The order topology induced by $<_\al$ coincides with the quotient topology. Clearly, each $q^\beta_\al$ is order preserving with respect to $<_\beta$ and $<_\al$.
This completes the proof.
\end{pf}

\begin{tw} There exists a $0$-dimensional linearly ordered Valdivia compact space $K$ with the following properties:
\begin{enumerate}
	\item[(1)] $K$ is nowhere first countable; in fact every nonempty open interval contains both a copy of $\omega_1+1$ and a copy of $1+\omega_1^{-1}$.
	\item[(2)] Every nonempty clopen subset of $K$ is order isomorphic to $K$.
	\item[(3)] Every linearly ordered Valdivia compact is an order preserving image of $K$.
\end{enumerate}
\end{tw}

\begin{pf} Let $C=\pair{2^\omega}{<_{\lex}}$ be the Cantor set, endowed with the lexicographic order. Denote by $Q$ the set of all $p\in C$ which are isolated either from the left or from the right (in particular $0,1\in Q$). Define $D=(C\times\sn0)\cup(Q\times C)$, endowed with the lexicographic order, and let $\map\pi DC$ be the projection. Observe that $D$ is order isomorphic to $C$ and $\pi$ is an order preserving retraction. A possible right inverse $\map iCD$ to $\pi$ can be defined by setting $i(p)=\min\pi^{-1}(p)$, in case where $p$ is isolated from the right, $i(p)=\max\pi^{-1}(p)$, in case where $p$ is isolated from the left and $i(p)=p$ if $p\in C\setminus Q$. 

Let $\S=\invsys C\pi{\omega_1}\xi\eta$ be the unique continuous inverse sequence satisfying $C_0=C$ and $\pi^{\xi+1}_\xi=\pi$ for every $\xi<\omega_1$, where each $C_\xi$ is identified with $C$ by using an order preserving isomorphism. The case of a limit stage $\delta<\omega_1$ makes no trouble, because $C_\delta$ is a $0$-dimensional, dense-in-itself linearly ordered metric compact, therefore order isomorphic to $C$. 

Define $K=\liminv\S$. Since $\pi$ is an order preserving retraction, $K$ is a linearly ordered Valdivia compact (by Theorem \ref{KM}). 

We show (2). Fix a nonempty clopen set $J\subs K$. Observe that $K+K$ is isomorphic to $K$, therefore we may assume that $J$ is a clopen interval $[a,b]$. We have that $J=\inv{\pi_\al}U$ for some clopen set $U\subs C_\al$. For convenience, we assume that $\al=0$. Since $\pi_0$ is order preserving, necessarily $U$ is an interval. Thus $U=[c,d]$, where $c=\pi_0(a)$ and $d=\pi_0(b)$. Let $J_\xi=\inv{(\pi^\xi_0)}U$ and let $p^\eta_\xi=\pi^\eta_\xi\rest J_\eta$. Then $\S_J=\invsys Jp{\omega_1}\xi\eta$ is a continuous inverse sequence with limit $J$. 
Finally, observe that $p^{\xi+1}_\xi$ is isomorphic to $\pi$. Indeed, if $x\in J_\xi$ is isolated from one of the sides in $J_\xi$ then it is also isolated from the same side in $C_\xi$. Further, $(p^{\xi+1}_\xi)^{-1}(x)=\pi^{-1}(x)$. Since $J_\xi$ is order isomorphic to $C_\xi$, this shows that $p^{\xi+1}_\xi$ is isomorphic to $\pi$.

Observe that (1) follows from (2) and (3), so it remains to show (3).
Fix a linearly ordered Valdivia compact $X$. By Proposition \ref{bnjerg}, we may assume that $\w(X)=\aleph_1$ and therefore $X=\liminv\invsys Xr{\omega_1}\xi\eta$, where the sequence is continuous, retractive and each $X_\xi$ is a metrizable $0$-dimensional compact space.
Since $X$ has another representation as the limit of a sequence of linearly ordered metric compacta with order preserving surjections, by the spectral theorem of \v S\v cepin (see \cite[Thm. 2]{Szczepin} or \cite[Section 2.2]{KM}), we may assume (replacing $\omega_1$ by a closed cofinal set) that each $X_\xi$ is linearly ordered and that each $r^\eta_\xi$ is order preserving. 

We construct inductively a sequence $\sett{f_\xi}{\xi<\omega_1}$ of order preserving surjections $\map{f_\xi}{C_\xi}{X_\xi}$ such that $f_\xi \pi^\eta_\xi = r^\eta_\xi f_\eta$ holds for every $\xi<\eta<\omega_1$. Clearly, the limit of such a sequence is the desired order preserving surjection from $K$ onto $X$.

We start with any order preserving surjection $\map{f_0}{C_0}{X_0}$. For a limit ordinal $\delta$, we define $f_\delta$ as the limit of $\setof{f_\xi}{\xi<\delta}$, using the continuity of both inverse sequences. Fix $\al<\omega_1$ and assume that $f_\xi$ have been constructed for $\xi\loe\al$. We construct $f_{\al+1}$ in such a way that $f_{\xi} \pi^{\al+1}_\xi = r^{\al+1}_\xi f_{\al+1}$ holds for $\xi=\al$. Then, using the inductive hypothesis, we deduce that this equality holds for every $\xi\loe\al$. Recall that $\pi^{\al+1}_\al=\pi$. Thus, up to order isomorphism, we have $C=C_\al$, $D=C_{\al+1}$. We set $Z:=X_\al$, $W:=X_{\al+1}$, $r:=r^{\al+1}_\al$ and $g:=f_\al$. In order to complete the proof, we need to find an order preserving surjection $\map hDW$ such that $rh = g\pi$, i.e. so that the following diagram commutes:
$$\xymatrix{
{D} \ar[r]^-{h}\ar[d]_{\pi} & {W}\ar[d]^{r}\\
{C} \ar[r]^-{g} & {Z}
}
$$
Given $z\in Z$ define $F_z=\inv{\pi}{g^{-1}(z)}$. Note that $F_z$ is a closed interval and hence, since $D$ is dense-in-itself, $F_z$ is either order isomorphic to the Cantor set or else $|F_z|\loe2$. Now observe that the required map $h$ must satisfy
\begin{equation}
\img h{F_z}=r^{-1}(z)\quad\text{ for every }z\in Z.
\tag{*}\end{equation}
Conversely, if $\map hDW$ satisfies (*) then it is an order preserving surjection such that $rh = g\pi$ holds. Thus, it suffices to show that $F_z$ has an order preserving map onto $r^{-1}(z)$ for every $z\in Z$. 

Fix $z\in Z$. Then $r^{-1}(z)$ is a closed metrizable interval in $W$ containing $z$ (recall that $Z\subs W$ and $r\rest Z=\id_Z$). 
If $r^{-1}(z)=\sn z$ then there is nothing to prove, so assume $r^{-1}(z)=[a,b]$, where $a<b$. Note that $Z\cap [a,b]=\sn z$.
Thus, if $z\ne b$ then $z$ is isolated from the right in $Z$ and consequently $\max g^{-1}(z)\in Q$. Similarly, if $z\ne a$ then $\min g^{-1}(z)\in Q$. In both cases, $g^{-1}(z)\cap Q\nnempty$ and hence $F_z=\inv\pi{g^{-1}(z)}$ is isomorphic to the Cantor set, therefore it has an order preserving map onto every linearly ordered metric compact, in particular onto $r^{-1}(z)$. This completes the proof.
\end{pf}

Let $L=\setof{x\in\Qyu^{\omega_1}}{|\suppt(x)|<\aleph_0}$ be endowed with the lexicographic ordering. The above space $K$
may also be described as the Stone space of the interval algebra over $L$. In fact, the set $L'\subs K$ consisting of all points which are isolated from the left-hand side is order isomorphic to $L$.
In contrast to Theorem \ref{wejqwfjqpw}, observe that the natural two-to-one order preserving quotient of $K$ produces a connected linearly ordered space in class \RC, which is nowhere separable. 

We do not know whether every linearly ordered continuous image of a Valdivia compact is an order preserving (or even just continuous) image of the space $K$ from the above theorem. It can be shown that every linearly ordered continuous image of a Valdivia compact space has weight $\loe\aleph_1$, by adapting the proof of Proposition \ref{bnjerg} and using the fact that $\omega_2+1$ is not an image of any Valdivia compact.
Below are other two open questions, which have been already mentioned in the text.

\begin{question}
Is class \R~stable under retractions?
\end{question}

\begin{question}
Is every space from $\RC$ a continuous image of a space from class \R?
\end{question}

\vspace{3mm}

{\bf Acknowledgments.} The author would like to thank the following institutions for their hospitality and support, which helped him very much to complete this work:
University of Prince Edward Island (Charlottetown, Canada);
The Fields Institute for Research in Mathematical Sciences (Toronto);
CNRS and Universit\'e Paris VII (France); 
Centre de Recerca Matem\`atica in Barcelona (Spain); Universitat de Val\`encia (Spain).
The author is grateful to the anonymous referee for pointing out several corrections and improvements.


\begin{thebibliography}{99}

\bibitem{A}
{\sc A. V. Arkhangel'ski\u\i},
{\em Topological function spaces\/},
Translated from the Russian by R. A. M. Hoksbergen. Mathematics and its Applications (Soviet Series), 78. 
Kluwer Academic Publishers Group, Dordrecht, 1992.

\bibitem{AMN} {\sc S. Argyros, S. Mercourakis, S. Negrepontis}, {\em Functional-analytic properties of Corson-compact spaces\/}, 
Studia Math. {\bf89} (1988), no. 3, 197--229.

\bibitem{Ba91a} {\sc I. Bandlow},
{\em A construction in set-theoretic topology by means of elementary substructures\/},
Z. Math. Logik Grundlag. Math. {\bf 37} (1991), no. 5, 467--480.

\bibitem{Ba91} {\sc I. Bandlow}, {\em A characterization of Corson-compact spaces\/}, Comment. Math. Univ. Carolinae {\bf 32}, 3 (1991) 545--550.

\bibitem{Ba92} {\sc I. Bandlow}, {\em On the absoluteness of openly-generated and Dugundji spaces\/}, Acta Univ. Carolin. Math. Phys. {\bf 33} (1992), no. 2, 15--26.

\bibitem{Ba94} {\sc I. Bandlow}, {\em On function spaces of Corson-compact spaces\/}, Comment. Math. Univ. Carolin. {\bf 35} (1994), no. 2, 347--356.

\bibitem{BKT} {\sc M. Burke, W. Kubi\'s, S. Todor\v{c}evi\'c}, {\em Kadec norms on spaces of continuous functions\/}, submitted
\footnote{Preprint available at \texttt{http://arxiv.org/abs/math.FA/0312013}.}.

\bibitem{DG} {\sc R. Deville, G. Godefroy}, {\em Some applications of projective 
resolutions of identity\/}, Proc. London Math. Soc. {\bf(3) 67} (1993), 
no. 1, 183--199.

\bibitem{DGZ} {\sc R. Deville, G. Godefroy, V. Zizler}, {\em Smoothness and Renormings in Banach Spaces\/}. Pitman Monographs and Surveys in Pure and Applied Mathematics, 64. Longman Scientific \& Technical, Harlow; copublished in the United States with John Wiley \& Sons, Inc., New York, 1993.

\bibitem{Dow} {\sc A. Dow}, {\em An introduction to applications of elementary submodels to topology\/}, Topology Proc. {\bf 13} (1988), no. 1, 17--72.

\bibitem{Kalenda_segment} {\sc O. Kalenda},
{\em Embedding of the ordinal segment $[0,\omega\sb 1]$ into continuous images of Valdivia compacta\/},
Comment. Math. Univ. Carolin. {\bf 40} (1999), no. 4, 777--783.

\bibitem{Kalenda} {\sc O. Kalenda}, {\em Valdivia compact spaces in topology and Banach space theory\/}, Extracta Math. {\bf 15} (2000), no. 1, 1--85.

\bibitem{KM} {\sc W. Kubi\'s, H. Michalewski}, {\em Small Valdivia compact spaces\/}, Topology Appl. (in press)\footnote{Preprint available at \texttt{http://arxiv.org/abs/math.GN/0507062}.}.

\bibitem{KU} {\sc W. Kubi\'s, V. Uspenskij}, {\em A compact group which is not Valdivia compact\/}, Proc. Amer. Math. Soc. {\bf 133} (2005), No. 8, 2483-2487.

\bibitem{Szczepin} {\sc E. V. \v S\v cepin}, {\em
Topology of limit spaces with uncountable inverse spectra\/} (Russian),
Uspehi Mat. Nauk 31 (1976), no. 5 (191), 191--226. [English translation: Russian Math. Surveys 31 (1976), no. 5, 155--191.]

\bibitem{To_chains} {\sc S. Todor\v cevi\'c}, {\it
Chain-condition methods in topology\/}, Topology Appl. {\bf101} (2000),
no.~1, 45--82. 

\end{thebibliography}
\end{document}